\documentclass[12pt, leqno]{amsart}

\newlength{\mylength}
\usepackage{graphicx}
\usepackage{amssymb,amsmath,amsthm,amscd}
\usepackage{enumitem}
\usepackage{hyperref}
\usepackage{verbatim}

\numberwithin{equation}{section}
\allowdisplaybreaks[4]

\theoremstyle{plain}
\newtheorem{lemma}{Lemma}[section]

\newtheorem{definition}[lemma]{Definition}
\newcommand{\Def}{\begin{definition}}
\newcommand{\edf}{\end{definition}}
\theoremstyle{definition}

\newcommand{\Ztwo}{\mathbb {Z}_2}

\title[Corrections to Nipp's tables]{Minor corrections to Nipp's tables of quaternary quadratic forms}
\author[C.-h. Lee]{Chul-hee Lee}
\address{School of Mathematics, Korea Institute for Advanced Study \\ Seoul 130-722, Korea}
\email{chlee@kias.re.kr}
\date{\today}

\begin{document}

\begin{abstract}
We make some corrections to Nipp's tables of positive definite integral quaternary quadratic forms. They only affect $p$-adic densities and $p$-adic Jordan splittings in the Appendix of Nipp's book. There are no errors found in Nipp's list of quaternary forms for each genus.
\end{abstract}
\maketitle

\section{Introduction}

In 1991, Gordon Nipp published a book containing computer-generated tables of primitive positive definite integral quaternary quadratic forms \cite{MR1118842}. The tables list a representative of each equivalence class of such forms for each discriminant up to 500. The book has an accompanying 3.5-inch floppy disk, with more tables up to discriminant 1732. These tables are currently available online \footnote{\url{http://www.math.rwth-aachen.de/~Gabriele.Nebe/LATTICES/nipp.html}}. 

The book has been a valuable resource for researchers for decades. It has been particularly useful for classifying low-rank quadratic forms with specific properties, allowing for exhaustive enumeration. It seems that there is a consensus on the completeness and accuracy of Nipp's tables, as they have been used by many researchers over a long period of time.

The book consists of three parts. The first is an introductory section that explains the method used to generate the tables. The second section has a printed version of the tables through discriminant 500, which is the main body of the book. In these tables, the forms are organized into genera, and for each genus, its mass is given. For an individual form belonging to a genus, its level, the Hasse symbol and the number of automorphisms are included. In the Appendix, there are tables of $p$-adic densities and $p$-adic Jordan splittings for each genus for $p=2$ and for each odd prime $p$ dividing the discriminant.

We recently made some attempts to write a computer code to find the canonical form of a quadratic form over $\Ztwo$, in the sense of \cite{MR0476636}. It was natural that we decided to use Nipp's table to test the accuracy of our implementation. Since the canonical forms of two equivalent forms over $\Ztwo$ are supposed to be the same, the canonical forms of a quadratic form from Nipp's table and its $2$-adic Jordan splitting in the Appendix should be identical. 

After we started testing all items in Nipp's tables, we could observe the complete agreement of canonical forms up to discriminant 1213. Quite interestingly, computer started to reveal some problematic items from discriminant 1216. After investigating possible sources of discrepancies, including errors in our computer code, we came to the conclusion that there are some genera with incorrect $p$-adic densities and $p$-adic Jordan splittings in the Appendix of \cite{MR1118842} and its accompanying disk.

Since they could have been caused by certain errors in Nipp's computer code for $p$-adic Jordan splittings, which could seriously affect the integrity of his tables, we have reexamined whether the quadratic forms in each genus do belong to the same genus. We also computed the mass of each genus by computing its $p$-adic densities. Fortunately, there are no errors found in Nipp's main tables. It was only the Appendix that contained incorrect information. Thus, there is no change in the classification by Nipp.

\section{List of genera with incorrect information}
Let us consider the first quadratic form of discriminant 1216 belonging to the genus with id \#15 in Nipp's table. Its coefficients are $$[f_{11},f_{22},f_{33},f_{44},f_{12},f_{13},f_{23},f_{14},f_{24},f_{34}]=[1, 1, 11, 11, 1, 0, 0, 1, 0, 8],$$ and the corresponding symmetric half-integral matrix is 
$$
\left(
\begin{array}{cccc}
 1 & 1/2 & 0 & 1/2 \\
 1/2 & 1 & 0 & 0 \\
 0 & 0 & 11 & 4 \\
 1/2 & 0 & 4 & 11 \\
\end{array}
\right).
$$
In the Appendix, its 2-adic Jordan splitting is given as $[2A]+[(58/3)+(38/29)]$, which is, in matrix form,
$$
\left(
\begin{array}{cccc}
 2 & 1 & 0 & 0 \\
 1 & 2 & 0 & 0 \\
 0 & 0 & 58/3 & 0 \\
 0 & 0 & 0 & 38/29 \\
\end{array}
\right).
$$
% Here, $A=\left(
% \begin{array}{cc}
%  1 & 1/2 \\
%  1/2 & 1 \\
% \end{array}
% \right)$.
However, these two matrices are not equivalent over $\Ztwo$. Its $2$-adic density is incorrectly given as $98304$, which should have been $3072$. The $p$-adic splitting for $p=19$, which divides 1216, is also incorrectly given.

The following table %Table \ref{table}
gives the complete list of genera with incorrect $p$-adic density, or incorrect $p$-adic splitting from Nipp's tables. For example, the genus in the above example is denoted by 1216\#15.
%\begin{tiny}
\begin{table}[ht]
\begin{tabular}{c|c|c|c|c|c|c}
% $(u_1,u_2,u_3)$ & $(u',K)$ & $\text{genus symbol}$ & $(u_1,u_2,u_3)$ & $(u',K)$ & $\text{genus symbol}$ \\
%\hline
1216\#15 &1216\#16 &1216\#17 &1216\#18 &1216\#19 &1216\#20 &1232\#16 \\
1280\#24 &1280\#25 &1280\#26 &1280\#27 &1280\#28 &1280\#29 &1280\#30 \\
1280\#31 &1280\#32 &1280\#33 &1280\#34 &1280\#35 &1280\#36 &1280\#37 \\
1280\#38 &1280\#39 &1280\#40 &1280\#41 &1280\#42 &1280\#43 &1296\#30 \\
1296\#31 &1296\#32 &1296\#33 &1296\#34 &1296\#35 &1296\#36 &1296\#37 \\
1296\#38 &1296\#39 &1296\#40 &1344\#26 &1344\#27 &1344\#28 &1344\#29 \\
1344\#30 &1344\#31 &1344\#32 &1344\#33 &1344\#34 &1344\#35 &1344\#36 \\
1344\#37 &1344\#38 &1360\#16 &1408\#24 &1408\#25 &1408\#26 &1408\#27 \\
1408\#28 &1472\#16 &1472\#17 &1472\#18 &1472\#19 &1472\#20 &1488\#16 \\
1536\#32 &1536\#33 &1536\#34 &1536\#35 &1536\#36 &1536\#37 &1536\#38 \\
1536\#39 &1536\#40 &1536\#41 &1536\#42 &1536\#43 &1536\#44 &1536\#45 \\
1536\#46 &1536\#47 &1536\#48 &1536\#49 &1536\#50 &1536\#51 &1536\#52 \\
1536\#53 &1536\#54 &1536\#55 &1536\#56 &1536\#57 &1536\#58 &1600\#25 \\
1600\#26 &1600\#27 &1600\#28 &1600\#29 &1600\#30 &1600\#31 &1616\#7 \\
1616\#8 &1620\#23 &1620\#24 &1620\#25 &1620\#26 &1620\#27 &1620\#28 \\
1620\#29 &1620\#30 &1620\#31 &1620\#32 &1620\#33 &1620\#34 &1620\#35 \\
1620\#36 &1664\#24 &1664\#25 &1664\#26 &1664\#27 &1664\#28 &1680\#28 \\
1680\#29 &1680\#30 &1680\#31 &1680\#32 &1728\#43 &1728\#44 &1728\#45 \\
1728\#46 &1728\#47 &1728\#48 &1728\#49 &1728\#50 &1728\#51 &1728\#52 \\
1728\#53 &1728\#54 &1728\#55 &1728\#56 &1728\#57 &1728\#58 &1728\#59 \\
1728\#60 &1728\#61 &1728\#62 &1728\#63 &1728\#64 &1728\#65 &1728\#66 \\
1728\#67 &1728\#68 &1728\#69 &1728\#70 &1728\#71 &1728\#72 &1728\#73 \\
1728\#74 &1728\#75 &1728\#76 &1728\#77 &1728\#78 &1728\#79 &1728\#80
\end{tabular}
\label{table}
\caption{List of genera with incorrect $p$-adic density or $p$-adic splitting}
\end{table}
%\end{tiny}
The corrected version of $p$-adic densities and $p$-adic splittings for these genera is available online at \url{https://github.com/chlee-0/nipp}.

% \section*{Acknowledgements}
% The author thanks ...

\bibliographystyle{amsalpha}
\bibliography{main}
\end{document}